\documentclass{amsart}
\usepackage[utf8]{inputenc}
\usepackage{todonotes  }
\usepackage{amsthm}
\usepackage{amsmath}
\usepackage{amssymb}
\usepackage{graphicx}
\usepackage{caption}
\usepackage{subcaption}
\usepackage{float}
\usepackage{dsfont}
\usepackage{longtable}
\newtheorem{theorem}{Theorem}[section]

\theoremstyle{definition}

\newcounter{x}
\usepackage{indentfirst}
\usepackage{mathalfa}
\usepackage{hyperref}

\newcommand{\st}{s.t.}
\newcommand\norm[1]{\left\lVert#1\right\rVert}
\newcommand{\Sp}{\textup{Sp}}
\newcommand{\py}{\textit{Python}}
\newcommand{\code}[1]{\texttt{#1}}

\newcommand\cM{\mathcal M}
\newcommand\cL{\mathcal L}
\newcommand\cT{\mathcal T}
\newcommand\cE{\mathcal E}
\newcommand\Id{\textup{Id}}

\newcommand{\sn}{\textup{sn}}
\newcommand{\cn}{\textup{cn}}
\newcommand{\amp}{\varphi}

\newcommand{\ie}{i.e.\ }
\newcommand{\eps}{\varepsilon}
\newcommand{\resp}{resp.\ }
\newcommand{\domains}{\mathcal{D}}
\newcommand{\modsep}{\,|\,}
\newcommand{\ifrac}[2]{#1/#2}
\newcommand{\mmc}{\varkappa}

\newtheorem*{conjecture}{Conjecture}
\newtheorem{proposition}[theorem]{Proposition}
\newtheorem{rmk}[theorem]{Remark}

\title[Spectral rigidity of ellipses]{Numerical evidence of Dynamical spectral rigidity of ellipses among smooth $\mathbb{Z}_2$-symmetric domains}
\author{Shanza Ayub and Jacopo De Simoi}
\date{\today}
\begin{document}
\maketitle
\begin{abstract}
  We present numerical evidence for spectral rigidity among
  $\mathbb{Z}_{2}$-symmetric domains of ellipses of eccentricity smaller than $0.30$.
\end{abstract}
\section{Introduction}
The famous question “Can one hear the shape of a drum?“ posed by
M. Kac in~\cite{csd} has motivated over 50 years of research into what
is now called the \emph{Inverse Spectral Problem}.  In this paper we
present numerical evidence to support a conjecture that is closely
related to this problem.

Let us introduce the main concepts so that we can present our results.
In this paper a \emph{domain} will refer to a subset
$\Omega\subset\mathbb{R}^{2}$ which is open, connected, bounded and
whose boundary $\partial\Omega$ is a sufficiently smooth curve; for
simplicity\footnote{ Our discussion can be actually applied to domains
  whose boundary is $C^{8}$-smooth, but we do not want to insist on
  this point here.} we consider here domains with $C^{\infty}$
boundary.  We denote with $\domains$ the set of all such domains.

Given a domain $\Omega\in\domains$, we denote its \emph{Laplace
  Spectrum} with
\begin{align*}
  \Sp(\Omega)\ &= \{0 < \lambda_0 \leq \lambda_1 \leq \cdots \leq \lambda_k \leq \cdots\}
\end{align*}
where the $\lambda_i$ are the eigenvalues of the Dirichlet\footnote{
  Historically, the majority results in the field have been obtained
  with Dirichlet boundary conditions, although other type of boundary
  conditions can be treated and are equally relevant.  In this paper
  we will follow this long established tradition and consider only
  Dirichlet boundary conditions. } Boundary Problem, \ie those $\lambda$
for which there exists $u\in L^{2}(\Omega)$ so that:
\begin{align*}
  \Delta u(x) + \lambda u(x) &= 0 \text{ if } x\in\Omega \\
  u(x)&= 0 \text{ if } x\in\partial\Omega.
\end{align*}

Kac's question can be then expressed, more formally, as “Does
$\Sp(\Omega)$ determine $\Omega$?”.  Clearly, domains that are
isometric to each other (\ie can be obtained by one another via a
composition of rotations and translations) will have the same Laplace
spectrum.  From now on we will, in this paper, identify isometric
domains, \ie we consider two domains to be equal if they are
isometric.  Two domains $\Omega$ and $\Omega'$ are said to be
\emph{Laplace isospectral} if $\Sp(\Omega) = \Sp(\Omega')$.  We can
thus further rephrase Kac's question as: “Are isospectral domains
necessarily isometric?”\\

In full generality, this question has a negative answer: in~\cite{gww}
the authors construct an explicit example of a pair of isospectral
domains that are not isometric, and many more domains can indeed be
constructed by similar methods.  However all such constructions yield
domains that are neither smooth nor convex.  In fact Kac's question is
still open if we require $\Omega$ to have a smooth boundary: this
problem is indeed notoriously hard.  In order to obtain some results
with the current technology, some further restrictions on the class of
admissible domains are needed.

In order to proceed with our discussion, let us introduce some further
notation.  Let us fix a class $\cM\subset\domains$ of domains; a
domain $\Omega\in\cM$ is said to be $\cM$-\emph{spectrally determined}
if there are no other domains in $\cM$ that are isospectral to
$\Omega$.  In other words, let us define the $\cM$-\emph{isospectral
  set of $\Omega$}:
\begin{align*}
  \textup{Iso}_{\cM}(\Omega) = \{\Omega'\in\cM \text{ so that } \Sp(\Omega) = \Sp(\Omega')\};
\end{align*}
then $\Omega$ is $\cM$-spectrally determined if
$\textup{Iso}_{\cM}(\Omega) = \{\Omega\}$.  Moreover, we say that
\emph{one can solve the Inverse Spectral Problem in $\cM$} if every
domain $\Omega \in \cM$ is $\cM$-spectrally determined.\\

A surprising relation exists between $\Sp(\Omega)$ and a dynamically
determined object, \textit{the Length Spectrum} of $\Omega$.  Let us
recall that, given a domain $\Omega$, one can consider the billiard
dynamics inside $\Omega$; some trajectories of this dynamics might be
\emph{periodic}.  Geometrically, periodic billiard trajectories
correspond to closed polygons, not necessarily convex, inscribed in
$\Omega$ with the property that, at each vertex, the angles that
either of the two sides that join at the vertex form with the tangent
line to the boundary are equal to each other (this is the well-known
law of optical reflection).  The perimeter of such a polygon is said
to be the \emph{length} of the corresponding periodic trajectory.  The
number of sides of such a polygon is called the \emph{period} of the
trajectory.  The Length Spectrum of $\Omega$ is then defined as the
set:
\begin{align*}
  \mathcal{L}(\Omega) = \mathbb{N}\{\text{length of all
  periodic billiard trajectories of $\Omega$}\} \cup
  \mathbb{N}\{|{\partial\Omega}|\}
\end{align*}
where $|{\partial\Omega}|$ is the length of the boundary of $\Omega$
and the factor of $\mathbb{N}$ accounts for the fact that one can
always consider a periodic orbit traversed multiple times as a
periodic orbit of length equal to a multiple of the original length of
the orbit.

The relation between the Laplace Spectrum and the Length Spectrum can
be stated as follows: consider the following distribution
 \begin{equation} \label{eq1}
w(t) := \sum_{\lambda_i \in \Sp(\Omega)} \cos{(t\sqrt{\lambda_i})}
 \end{equation}
 which is called the \textit{wave trace distribution}; then it has
 been proved in~\cite{pgr} that:
 \begin{equation} \label{wtr}
\text{sing supp }(w(t)) \subseteq \pm \mathcal{L}(\Omega) \cup \{0\}
\end{equation}
where $\text{sing supp }(w(t))$ denotes the \textit{singular support}
of $w(t)$.  Moreover, if $\Omega$ satisfies some generic conditions,
which can be expressed purely in dynamical terms, it has been shown
that the inclusion in~\eqref{wtr} is indeed an equality: in
particular, for generic domains, the Laplace Spectrum determines the
Length Spectrum (see e.g.~\cite[Remark 2.10]{dsr} and references
therein).  And so, just how the Inverse Spectral Problem has been
posed in terms of the Laplace Spectrum, one can set up an Inverse
Problem in terms of the Length Spectrum.  A domain
$\Omega \in \mathcal{M}$ is said to be $\cM$-\textit{dynamically
  spectrally determined} if $\Omega$ is the unique element of
$\mathcal{M}$ with the same Length Spectrum.  Hence, the Inverse
Dynamical Problem asks: is every $\Omega \in \mathcal{M}$ dynamically
spectrally determined?\\

The inverse problem (either in the Laplace formulation or the
dynamical formulation) turns out to be extremely hard (see the next
section for some available prior results).  A related question, that
proved to be more tractable, is the problem of \emph{spectral
  rigidity}: we say that a domain $\Omega$ is $\cM$-spectrally rigid
(\resp $\cM$-dynamically spectrally rigid) if every $C^{1}$-family of
domains $\{\Omega_{t}\}_{t\in(-\eps,\eps)}$ with $\Omega_{0} = \Omega$
and with the property that $\Sp(\Omega) = \Sp(\Omega_{t})$ (\resp
$\cL(\Omega) = \cL(\Omega_{t})$) for any $t\in(-\eps,\eps)$ is a
trivial family (\ie a family of isometric domains).  In other words, a
domain $\Omega$ is $\cM$-spectrally rigid (\resp $\cM$-dynamically
spectrally rigid) if every $C^{1}$-deformation in $\cM$ preserving the
Laplace (\resp Length) Spectrum is necessarily a trivial deformation.
Clearly, if a domain is $\cM$-spectrally determined, then it is
$\cM$-spectrally rigid, but the converse, in general, is not guarantee
d to be true.

In this paper we provide some numerical evidence to support the following
conjecture:
\begin{conjecture}
  Ellipses are dynamically spectrally rigid among axially-symmetric
  smooth convex domains.
\end{conjecture}
We will present our results in full detail in
Section~\ref{sec:results-concl-remark}.

\subsection{Related prior results.} 

It has been shown in~\cite{zelditch} that the inverse spectral problem
can be solved in a class $\cM$ of domains that are convex,
axially-symmetric, analytic and satisfy a generic dynamical condition.
However, such results depend crucially on the analyticity assumption,
and cannot be extended in any way to the case of smooth domains.  Some
progress was made in~\cite{csr}, where Hezari and Zelditch showed that
if $\Omega_0$ is an ellipse and $\Omega_\tau$ is a $C^1$ Dirichlet (or
Neumann) isospectral deformation of $\Omega_0$ through $C^\infty$
domains which preserves the $\mathbb{Z}_2 \times \mathbb{Z}_2$
symmetry group of the ellipse, then it is necessarily \textit{flat}
(i.e. all derivatives must vanish for $\tau= 0$). This result shows
that ellipses are \textit{infinitesimally spectrally rigid} among
$C^\infty$ domains with the symmetries of the ellipse.  Very recently,
in~\cite{hz2}, the same authors proved that ellipses of small
eccentricity are spectrally determined among all $C^{\infty}$ domains.
This, of course, settles the conjecture that we are investigating for
small values of the eccentricity.  However, no bound on the
\emph{smallness} of eccentricity is provided.  We point out that all
studies mentioned thus far have used more traditional (\ie non
dynamical) approach to Laplace spectral rigidity and determination
(except in~\cite{hz2}, where dynamical results obtained in~\cite{adk}
and~\cite{bcc} are crucially employed).


In this paper we rely on the dynamical technique used in~\cite{dsr} to
investigate numerically the problem.  In~\cite{dsr}, the authors prove
the \textit{dynamical spectral rigidity} of $\mathbb{Z}_2$-symmetric
strictly convex domains close to a circle.  The proof hinges on the
construction, for each $\mathbb{Z}_2$-symmetric strictly convex domain
$\Omega$, of an operator called the \textit{linearized isospectral
  operator}, which we denote with $\mathcal{T}_{\Omega}$.
Then they proved that the injectivity of $\mathcal{T}_{\Omega}$
implies the dynamical spectral rigidity of $\Omega$ among
$\mathbb{Z}_2$-symmetric strictly convex domains; finally they prove
that if $\Omega$ is sufficiently close to a circle, then the operator
$\mathcal{T}_{\Omega}$ is injective.

In this paper, we use the method outlined above and compute,
numerically, the linearized isospectral operator for ellipses of
various eccentricities; we then check numerically whether or not this
operator is injective to give evidence for the spectral rigidity of
the corresponding ellipse.

\textbf{Outline of the paper.} In
Section~\ref{sec:ellips-ellipt-integr}, we recall, for the readers'
convenience, some definitions regarding elliptic integrals and
elliptic functions.  In Section~\ref{sec:bill-dynam-ellips}, we
describe the method discussed in~\cite{dsr} in more detail and explain
how it will be used for the purposes of this paper.  In
Section~\ref{sec:results-concl-remark}, we present and discuss our
numerical results as based on the method described in
Section~\ref{sec:bill-dynam-ellips}.  We finally record in the tables
in the appendix the values obtained by our numerical investigation.

\section{Ellipses and elliptic integrals}\label{sec:ellips-ellipt-integr}
In this section, we recall a few basic definitions and set some
important notation used in the rest of the paper.  An ellipse centered
at the origin with semi-axes oriented along the coordinate axes and of
lengths $0<b\leq a$ is defined as follows:
\begin{align*}
\mathcal{E}_{a,b} = \Bigg \{(x,y) \in \mathbb{R}^{2} : \frac{x^2}{a^2}
  + \frac{y^2}{b^2} = 1 \Bigg \}.
\end{align*}
The value $a$ is the length of the major semi-axis and $b$ is the
length of the minor semi-axis. We denote the eccentricity of the
ellipse with $e = \sqrt{1-\ifrac{b^2}{a^2}} \in [0,1)$; it is a
measure of how close the ellipse is to a circle\footnote{Note that an
  ellipse with $e=0$ is a circle.}.  We let
$c = a\cdot e = \sqrt{a^2 - b^2}$ denote the distance from the center
to any of the foci.

We further recall the definitions for Elliptic Integrals and Jacobi
Elliptic functions (we also refer the reader to~\cite{bcc} and
references therein for a more comprehensive presentation). For $m\in[0,1)$:
\begin{itemize}
\item \textit{Incomplete elliptic integral of the first kind}: for
  $\varphi\in [0,\pi/2]$ we let
      \begin{align*}
        F(\varphi\modsep m) := \int_0^\varphi \frac{1}{\sqrt{1-m\sin^2{\varphi'}}} d\varphi';
      \end{align*}%
      here $k = \sqrt{m}$ is called the \textit{modulus}; the quantity
      $\varphi$ is called the \textit{amplitude}.
    \item \textit{Complete elliptic integral of the first
        kind}:
      \begin{align*}
        K(m) = F({\pi}/{2}\modsep m).
      \end{align*}
    \item \textit{Incomplete elliptic integral of the second kind}:
      \begin{align*}
        E(\varphi\modsep m) &= \int_0^{\varphi}\sqrt{1-m\sin^{2}(\varphi')}d\varphi'.
      \end{align*}
    \item \textit{Complete elliptic integral of the second kind}:
      \begin{align*}
        E(m) = E(\pi/2\modsep m).
      \end{align*}
\end{itemize}

We also recall the \textit{Jacobi elliptic functions} $\sn(u,m)$ and
$\cn(u,m)$ which are obtained by inverting the incomplete elliptic
integrals of the first kind, where $u$ is called the \textit{argument}
and $m$ is the \textit{modulus} as described earlier.  They are
defined in such a way that for a given $m$, we have
$\cn(u,m) = \cos(\varphi)$ and $\sn(u,m) = \sin(\varphi)$ where
$\varphi$ is so that $F(\varphi\modsep m) = u$; $\varphi$ is called
the \emph{amplitude} of $u$ and will be denoted by $\amp(u\modsep m)$.

\section{Billiard Dynamics in ellipses and Methods}\label{sec:bill-dynam-ellips}

In this section we provide definitions relating to the billiard map
within an ellipse $\mathcal{E}$, we describe the method described
in~\cite{dsr} for determining spectral rigidity in our study and we
present the details of the implementation of this method.

More precisely, in Section~\ref{sec:method}, we will describe the
strategy, proposed in~\cite{dsr}, that we will implement for
determining spectral rigidity among $\mathbb Z_{2}$-symmetric domains;
then, in Section~\ref{sec:ellipses}, we will specialize our discussion
to the case of ellipses; finally, in Section~\ref{sec:numerical}, we
describe how the numerical computations were carried out in this
study.

\subsection{The strategy}\label{sec:method}
Let $\Omega$ be a $\mathbb Z_{2}$-symmetric domain; to fix ideas we
assume that the perimeter $|\partial \Omega| = 1$; since domains are
defined up to rigid motions, we can assume the symmetry axis of
$\Omega$ to coincide with the $x$-axis.  By convexity,
$\partial \Omega$ intersects the symmetry axis in two points: let us
choose one of them (arbitrarily) and denote it by $P$.  Let us further
assume (by possibly translating $\Omega$ along the $x$-axis) that the
other intersection point of $\partial\Omega$ with the $x$-axis is $-P$.

\begin{proposition}[{see~\cite[Lemma 4.3]{dsr}}]\label{p_existence-orbits}
  For any $q > 2$, there exist $q$ points $(X_{0}^{q},\cdots, X_{q-1}^{q})$,
  with $X_{j}^{q}\in\partial\Omega$, $X_{0}^{q} = P$, which correspond
  to collision points of a periodic billiard orbit of period $q$,
  with the property that the polygonal curve $X_{0}^{q}X_{1}^{q}\cdots
  X_{q-1}^{q}$ is a simple closed curve.
\end{proposition}
\begin{rmk}\label{r_bouncing-ball}
  If $q = 2$, due to symmetry and smoothness of $\partial\Omega$, we
  are guaranteed that $\partial\Omega$ intersects the $x$-axis at a
  right angle; therefore taking $X_{0}^{2} = P$ and $X_{1}^{2} = -P$
  yields, in fact, an orbit of the billiard table, which is called the
  \emph{bouncing ball} orbit.  This orbit is, strictly speaking, not
  \emph{simple}, but will be considered as such in the sequel.  The
  orbits of larger period can always be constructed by variational
  methods (see e.g.~\cite[Lemma 4.3]{dsr}).
\end{rmk}
In order to be consistent, if several orbits of such type exist, we
choose the points $X_{j}^{q}$ so that they give an orbit of maximal
length and so that the length of $X^{q}_{0}X^{q}_{1}$ is maximal among
the orbits of maximal lengths.  Moreover, for $q = 1$, it turns out to
be convenient to define, conventionally, $X_{0}^{1} = P$.

Using this sequence of periodic orbits, we proceed to construct the
associated so-called \emph{linearized isospectral operator}.  In order
to do so we need to introduce some more notation.  First, by
definition of billiard, each collision point $X_{j}^{q}$ is such that
the angle between the incoming edge $X_{j-1}^{q}X_{j}^{q}$ and the
tangent vector to $\partial\Omega$ at $X_{j}^{q}$ equals the angle
between the outgoing edge\footnote{ The subscripts in $X^{q}_{j}$ are
  considered to be modulo $q$.} $X_{j}^{q}X_{j+1}^{q}$ and the same
tangent vector.  Let us denote this angle by $\phi_{j}^{q}\in(0,\pi)$.
Conventionally, for $q = 1$, we define $\phi_{0}^{1} = \pi/2$.  Next,
we introduce a convenient parametrization of $\partial\Omega$ as
follows.  Let $s$ denote the arc-length parametrization of
$\partial\Omega$ with the choice of origin so that $P$ corresponds to
$s = 0$ and let $\rho(s)$ denote the \textit{radius of curvature} of
$\partial\Omega$ expressed {in terms of $s$.}  Then we define the
\textit{Lazutkin parameterization}, denoted by $x$, as follows:
\begin{align*}
x(s) &= C\int_0^s\rho(s)^{-2/3} ds
&\textrm{where } C &= \Bigg
  [\int_{0}^{1}\rho(s)^{-2/3}ds \Bigg ]^{-1}.
\end{align*}
  We also define the \textit{Lazutkin weight}
function:
\begin{align*}
  \mu(x) = \frac{1}{2C\rho(x)^{1/3}},
\end{align*}
where $\rho(x)$ denotes the radius of curvature as a function of the
Lazutkin coordinate $x$.  We refer the reader to \cite[Appendix
A.2]{dsr} for an in-depth discussion of properties of this
parameterization.  We then denote with $x_{j}^{q}$ the Lazutkin
coordinate of the collision point $X_{j}^{q}$ (in particular,
$x_{0}^{q} = 0$ for any $q\ge1$).

Then, for any $q\ge 1$ and $j\ge 1$ we define the following
quantity:
\begin{align}\label{eq:Tqj-formula}
\mathcal{T}_{q,j} &= \sum_{n=0}^{q-1}\frac{\cos(2\pi j
  x_n^q)}{\mu(x_n^q)} \sin(\phi_n^q)
\end{align}
\begin{rmk}
  The quantities $\cT_{q,j}$ admit the following geometrical
  interpretation.  Consider the deformation of domain $\Omega$ by the
  infinitesimal \emph{normal} perturbation described by the function
  $n(x)$; in other terms, at the point identified by $x$, we are
  deforming $\partial\Omega$ along the normal direction by $n(x)$
  (outward if $n(x) > 0$ or inward if $n(x) < 0$).  If $n$ preserves
  the $\mathbb Z_{2}$-symmetry (\ie it is an even function of $x$),
  the deformed domain will also be $\mathbb Z_{2}$-symmetric; assuming
  that the orbits found above persist this deformation, they will
  possibly change their length.  For $q \ge 2$ the quantity
  $\cT_{q,j}$ corresponds to the variation of the length of the $q$-th
  orbit by the deformation $n(x) = \cos (2\pi j x)$ (the $j$-th Fourier
  Mode).
\end{rmk}

The criterion for spectral rigidity of $\Omega$ that was proposed in
~\cite{dsr} can be now loosely stated as follows: if the infinite
matrix $\cT_{q,j}$ is not degenerate, then $\Omega$ is spectrally
rigid.  In order to properly state the non-degeneracy condition, we
find necessary to introduce some other notions.

Remarkably, the following proposition holds (see~\cite[Lemma B.1 and
Lemma 5.3]{dsr}):
\begin{proposition}
For any $j\ge 1$, the quantities
\begin{align}\label{eq:marvizi-melrose-def}
 \mmc_{j} = \lim_{q \to \infty} q^2 \mathcal{T}_{q,j}
\end{align}
exist and are finite.
\end{proposition}
\begin{rmk}
  In the notation of~\cite{dsr}, we have the expression
  \begin{align*}
    \mmc_{j} &= \tilde\ell_{\bullet}(\cos(2\pi jx)).
  \end{align*}
  The values of $\mmc_{j}$ are related to the variation of the first
  Marvizi--Melrose (see~\cite{scp}) coefficient of $\Omega$ by a
  perturbation that is the $j$-th Fourier harmonic in Lazutkin
  coordinates.
\end{rmk}

We then define the reduced matrix:
\begin{align*}
  \tilde\cT_{q,j} = \cT_{q,j}-\frac{\mmc_{j}}{q^{2}},
\end{align*}
and, for $\gamma > 0$, the sequence spaces:
\begin{align*}
  h_{\gamma} = \{b = (a_i)_{i\ge 0}\in\ell^\infty\st a_{0} =
  0\text{\ and }\lim_{j\to\infty}j^\gamma a_j = 0\}
\end{align*}
equipped with the norm $|b|_{\gamma} = \max_{j \ge0} j^\gamma|a_j|$.
Then (see~\cite[Lemma 5.3]{dsr}) the following proposition holds:
\begin{proposition}
  Let $\gamma\in(3,4)$; then the matrix $\tilde\cT$ acts as an
  operator $\tilde\cT:h_{\gamma}\to h_\gamma$.
\end{proposition}
Then the following criterion for spectral rigidity holds:
\begin{theorem}
  If $\tilde\cT$ is injective for some choice of $\gamma$, then
  $\Omega$ is spectrally rigid among $\mathbb Z_{2}$-symmetric
  domains.
\end{theorem}
The above is the criterion that we will use to determine whether or
not we can say that an ellipse of given eccentricity is spectrally
rigid.  We thus need to find an explicitly computable condition for
$\tilde\cT_{q,j}$ to be injective; we will use the following:
\begin{proposition}\label{p_injectivity}
  Let $X$ be a Banach space, $T:X\to X$ be an operator and $S:X\to X$
  to be an invertible operator; if $\|T-S\| < \|S^{-1}\|^{-1}$, where
  $\|\cdot\|$ denotes the operator norm, then $T$ is invertible; in
  particular $T$ is injective.
\end{proposition}
In our case, choosing $X = h_{\gamma}$, $T = \tilde\cT_{q,j}$ and
$S = \Id$, we can express the operator norm as:
\begin{align}
\norm{\tilde\cT-\Id} &= \max_{q}\Big(q^{\gamma}\sum_{j = 0}^{\infty}j^{-\gamma}|\tilde\cT_{q,j}-\Id_{q,j}|\Big)\notag\\\label{eq:norm-computaton}
&=
\max_{q}\Big(q^{\gamma}\sum_{j = 1}^{\infty}j^{-\gamma}\mid
\mathcal{T}_{q,j}-\delta_{qj}-\frac{\varkappa_j}{q^2} \mid\Big) \quad
\textrm{where} \quad \delta_{qj} = \left\{
\begin{array}{ll}
            1 & \quad j = q \\
            0 & \quad \textrm{otherwise}
\end{array}
                \right.
\end{align}
In the next section we specialize the above discussion to the
setting of elliptical domains.
\subsection{Dynamics in an elliptical billiard
  table}\label{sec:ellipses}
To fix ideas, we fix the family of ellipses in such a way that, for
each ellipse in the family, its center lies at the origin, its major
axis lies along the x-axis (observe that the major axis is a symmetry
axis for $\cE$) and its circumference (or perimeter), denoted by
$|\mathcal{E}|$, equals $1$.  We choose $P$ to be the
point $(a,0)\in\mathcal E$, that is the intersection of $\mathcal E$
with the positive $x$ semi-axis.

\begin{rmk}
Since ellipses have a $\mathbb Z_{2}\times\mathbb Z_{2}$ symmetry, we
could also consider the minor axis as an axis of symmetry.  This
choice seems not to affect our results. See Item~\ref{item:minor-axis}
in Section~\ref{sec:results-concl-remark}.
\end{rmk}

Given an ellipse, we consider its parametrization by the
\emph{amplitude} $\varphi$ as
\begin{align*}
  X(\varphi) = (a\sin(\varphi),-b\cos(\varphi)).
\end{align*}
Observe that with this choice of $\varphi$, the point $X(0)$
corresponds to the lowest point $(0,-b)$ of the ellipse.  In
particular\footnote{ This non-standard choice of origin for the
  parametrization $\varphi$ simplifies some formulae in Lazutkin
  coordinates.} we have $P = X(\pi/4)$.

The \emph{arc-length parametrization}, denoted by $s$, can then be
obtained observing that
$\frac{ds}{d\varphi} =
\sqrt{a^{2}\cos^{2}(\varphi)+b^{2}\sin^{2}(\varphi)} =
a\sqrt{1-e^{2}\sin^{2}(\varphi)}$; since $P$ corresponds to $s = 0$,
we conclude:
\begin{align*}
  s(\varphi)&= a (E(\varphi\modsep e^{2})-E(e^{2}));
\end{align*}
in particular we have:
\begin{align*}
  |\mathcal{E}| = 4aE(e^2).
\end{align*}
In particular, if we fix $e\in[0,1)$ and $|\mathcal{E}| = 1$, the
value of $a$ is determined from the above expression; $b$, could then
be found using the definition of $e$.  Our family of ellipses will
then parameterized by their eccentricity.

The first ingredient in the analysis recalled in
Section~\ref{sec:method} is the determination of a sequence of
periodic orbit of increasing period.  This is particularly convenient
to do inside an ellipse; we follow the approach in~\cite{bcc} and we
refer to their work for further details.  Here we limit ourselves to
state their observations without proofs.

In our discussion, we assume $e$ (and therefore $\mathcal E$) to be
fixed once and for all.

Let us start by considering the family of ellipses confocal to
$\mathcal E$ and contained within $\mathcal E$; such family of
ellipses can be parametrized by $ 0<\lambda < b$ as follows:
\begin{align*}
\mathcal{C}_\lambda = \Bigg \{ (x,y) \in \mathbb{R}:
\frac{x^2}{a^2-\lambda^2} + \frac{y^2}{b^2-\lambda^2} = 1 \Bigg \}.
 \end{align*}
Notice that for $\lambda\to 0$, $\mathcal C_{\lambda}$ approaches
$\mathcal E$ and for $\lambda\to b$, $\mathcal C_{\lambda}$ approaches
the segment joining the two foci.

Let us fix an orientation for $\mathcal C_{\lambda}$; to fix ideas we
will always take the counterclockwise orientation.  For any
$\lambda\in(0,b)$, let us proceed with the following inductive
geometrical construction: let $P_{0} = P$; then given
$P_{n}\in\mathcal E$, since $\mathcal C_{\lambda}$ is convex, there
exists a unique point $P_{n+1}\in\mathcal E$, so that the segment
$P_{n} P_{n+1}$ is tangent to $\mathcal C_{\lambda}$ and is
co-oriented with $C_{\lambda}$.  It turns out\footnote{ This fact is
  peculiar for elliptical billiards and follows from Poncelet's
  Porism} that the points $\{P_{n}\}$ form a sequence of collision
points of a billiard orbit.  The ellipse $\mathcal C_{\lambda}$ is a
\emph{caustic} for this billiard orbit.  Now, given $n > 0$, we can
define $p_{n}$ the \emph{winding number of $P_{n}$} as the number of
times that the polygonal path $P_{0}P_{1}P_{2}\cdots P_{n}P_{0}$ winds
around the origin\footnote{ A more precise definition of the winding
  number can be given, but we avoid giving it here, since this
  definition will suffice for our uses below.}.  For any given orbit
it can be proved (see for instance~\cite{stb}) that the ratio
$p_{n}/n$ converges to some number $\omega\in(0,1/2)$ that we call the
\emph{rotation number of the orbit}.  For some particular values of
$\lambda$, the corresponding orbit is \emph{periodic}; in particular
the rotation number is rational and it is given by the ratio of the
winding number of the orbit $p_{q}$ and its period $q$.  For the
purposes of our study, we will only consider periodic orbits with
rotation number $1/q$.  Such orbits are so that their trajectory is a
simple closed curve (recall Proposition~\ref{p_existence-orbits}).

Following~\cite{bcc}, let us define
$m_{\lambda} = \frac{a^2-b^2}{a^2-\lambda^2}$; notice that by
definition we have $m_{\lambda}\in(e^{2},1)$.  Then the rotation
number $\omega_{\lambda}$ of the orbit associated to the caustic
$\mathcal C_{\lambda}$ is (see~\cite{bcc}):
\begin{align}\label{eq:omega_lambda}
  \omega_\lambda = \frac{\delta_{\lambda}}{4 K(m_{\lambda})} \quad
  \textrm{where} \quad
  \delta_\lambda = 2F(\sin^{-1}(\lambda/b),m_{\lambda}).
\end{align}
Correspondingly, the collision points $P_{n}
= (x_{n},y_{n})$ can be obtained by the formula:
\begin{align}\label{eq:collision-points}
  X_{n} &= (a\cdot \sn(u_n,m_{\lambda}),-b\cdot \cn(u _n,m_{\lambda}))
        \quad\text{where $u_n = K(m_{\lambda})+n\delta_\lambda$.}
\end{align}
The above formula yields the collision points corresponding to the
given choice of $\lambda$.  In order to find the periodic orbits
described by Proposition~\ref{p_existence-orbits}, it thus suffices,
for any $q$, to find $\lambda_{q}$ so that
$\omega_{\lambda_{q}} = 1/q$.  Correspondingly
$\delta_{\lambda_{q}} = 4K(m_{\lambda_{q}})/q$ and in particular we obtain
\begin{align*}
  u_{n}^{q} = 4K(m_{\lambda_{q}})\left(\frac nq+\frac14\right).
\end{align*}

Observe that the periodic orbit of period $2$, corresponding to
rotation number $1/2$ cannot be obtained by this process, as it
corresponds to the limiting case $\lambda\to b$; however, it
corresponds to the orbit bouncing along the major axis (see
Remark~\ref{r_bouncing-ball}).  As in the previous section, we let
$X_{j}^{q}$ denote the points in the orbit of period $1/q$; from the
expression~\eqref{eq:collision-points} it is simple to obtain
(numerically) $\phi_n^q$.  We now need to find the value of the
Lazutkin parametrization $x_{j}^{q}$ at the collision points.  Such
coordinates could also be easily found numerically, but in the case of
ellipses, the Lazutkin parametrization can be analytically expressed
in terms of Elliptic Integrals: in fact, recall the expression for the
radius of curvature of an ellipse:
\begin{align*}
  \rho(\varphi) = \frac {a}{\sqrt{1-e^{2}}}(1-e^{2}\sin^{2}(\varphi))^{3/2}.
\end{align*}
Then substituting the above formula in the expresson for the Lazutkin
coordinate $x$, we conclude:
\begin{align*}
  x(\varphi)&= \frac{F(\varphi\modsep e^{2})-K(e^{2})}{4K(e^{2})}\quad
              \text{ and }C =
              \left[{4K(e^{2})a^{1/3}(1-e^{2})^{1/3}}\right ]^{-1}.
\end{align*}
and
\begin{align*}
  \mu({\varphi})&= 2 K(e^{2})\sqrt{\frac {{1-e^{2}}}{{1-e^{2}\sin^{2}\varphi}}}.
\end{align*}

In particular, for the periodic orbit of period $q$ we conclude:
\begin{align}\label{eq:collision-lazutkin}
  x_{n}^{q} = \frac
  {F(\varphi(u^{q}_{n}\modsep m_{\lambda_{q}}) \modsep e^{2})-K(e^{2})}
  {4K(e^{2})}
\end{align}
Observe that, as $q\to\infty$, $m_{\lambda_{q}}\to e^{2}$ and
$x_{n}^{q}\to n/q$.

\subsection{Details of our numerical
  implementation}\label{sec:numerical}
Here we list the tasks that we implemented numerically.  The code has
been written in \py{} and can be found in this
\href{https://github.com/shanzaayub/Spectral-Rigidity-Ellipse}{GitHub
  repository}.

First of all,
we considered orbits of period up to \code{maxq} and found that
setting \code{maxq}$ = 500$ was sufficient for our purposes.

\begin{itemize}
\item The first task that needs to be implemented is to find the
  sequence of values $\lambda_{q}$ corresponding to periodic orbits of
  rotation number $1/q$; we performed this computation by numerically
  inverting the formula~\eqref{eq:omega_lambda} for $\omega_{\lambda}$
  using the bisection method; this computation is implemented in the
  \py{} method \code{find\_lambda}.
\item Once the $\lambda_{q}$'s are found, we proceed to
  apply~\eqref{eq:collision-points} and find the collision points
  $X_{n}^{q}$; such points are stored as amplitudes $\varphi_{q}^{k}$
  and are calculated in the \py{} method
  \code{find\_collision\_pts}.
\item We need then, given the value $e$ of eccentricity, to compute
  the values $\cT_{q,j}$ for sufficiently many $q$'s and $j$'s
  according to~\eqref{eq:Tqj-formula}.  To this end we need to obtain
  the angles $\phi_{q}^{n}$; such computations are carried out by the
  \py{} method \code{sinphi\_lst}.  The method \code{T\_of\_q\_j}
  computes the values of the matrix $\cT_{q,j}$.
\end{itemize}

Once these methods are in place, we proceed with the actual
computation of the norms.  Here we report in detail our strategy:
\begin{itemize}
\item We fix an eccentricity $e\in(0,1)$; $e$ will sample various
  intervals in $(0,1)$, depending on the chosen value of $\gamma$.
\item Since the values $\mmc_{j}$ (defined
  in~\eqref{eq:marvizi-melrose-def}) do not depend on $q$, we cache
  them in the list \code{lambda\_MM}; observe that, due to the extra
  symmetries of the ellipse, $\mmc_{j} = 0$ if $j$ is odd.  In order
  to approximate the evaluation of the
  limit~\eqref{eq:marvizi-melrose-def}, we compute a term with
  sufficiently large $q$, until the difference of the values for
  $q$ and $q+1$ falls below a certain threshold (we took
  $10^{-6}$).
\item In order to compute numerically~\eqref{eq:norm-computaton},
we need to truncate the series to some value $J$, and stop at some $Q$
when computing the $\max$.  Below we explain how the two cutoffs $J$
and $Q$ have been chosen in our investigation:
\begin{itemize}
\item due to the factor $q^{\gamma}j^{-\gamma}$, the larger $q$ is,
  the higher the number of $j$'s that contribute substantially to the
  sum.  We decided to cut off the series at $J = C\cdot q$.  Observe
  that dropping the tail after this cutoff leads to an error of order
  $C^{-\gamma+1}$. In our computations we took $C = 100$.
\item the choice of the cutoff $Q$ depends on the eccentricity by
  means of the following argument: as $q$ increases, the values of
  $\cT_{q,j}$ should approach to the corresponding values for the case
  of a circle.  It is not difficult to compute such values explicitly:
  for the circle $\mmc_{j} = 0$ for any $j$ and $\cT_{q,j}$ can be
  computed explicitly (see~\cite{dsr}); let us denote the operator for
  the circle as $\cT_{q,j}^{0}$; we then obtain that
  \begin{align*}
    q^{\gamma}\sum_{j = 1}^{\infty}j^{-\gamma}\mid
\mathcal{T}^{0}_{q,j}-\delta_{qj}\mid &= |c_{q}-1| + c_{q}\sum_{s = 2}^{\infty}s^{-\gamma} =\\&= 1+c_{q}(\zeta(\gamma)-2),
  \end{align*}
  where $\zeta$ is the Riemann Zeta function and
  \begin{align*}
    c_{q} =
    \begin{cases}
      \pi^{-1}&\text{ if $q = 1$}\\
      \frac{\sin(\pi/q)}{\pi/q}&\text{ otherwise}.
    \end{cases}
  \end{align*}
  We then stop our computations as soon as the value of the $q$-th entry is
  within a reasonable accord with the value for the circle.  We fix
  this accord to be $10\%$.
\end{itemize}
\item the choice of the parameter $\gamma$ has been made to find a
  compromise with the computational time.  In fact, as $\gamma\to 3$,
  the quantity $\|\tilde\cT_{q,j}-\Id\|_{\gamma}$ tends to decrease,
  but computational times tend to increase. Hence, we fixed
  $\gamma = 3.5$ for a broad range of eccentricities and then computed
  the norms with $\gamma = 3.1$ and $\gamma = 3.01$ in narrower
  ranges.
\end{itemize}

\section{Results and Concluding Remarks}\label{sec:results-concl-remark}
We can now present the numerical evidence we collected by the methods
described in Sections 3.  This data suggests that the following
result should hold.
\begin{conjecture}
  Ellipses of eccentricity $e\in(0,0.3)$ are dynamically spectrally
  rigid among $\mathbb Z_{2}$-symmetric smooth convex domains.
\end{conjecture}
The value $0.3$ is certainly an artifact of our method and does not
represent a natural threshold.

\begin{figure}[ht]
\begin{subfigure}{0.5\textwidth}
\includegraphics[width=0.9\linewidth, height=5cm]{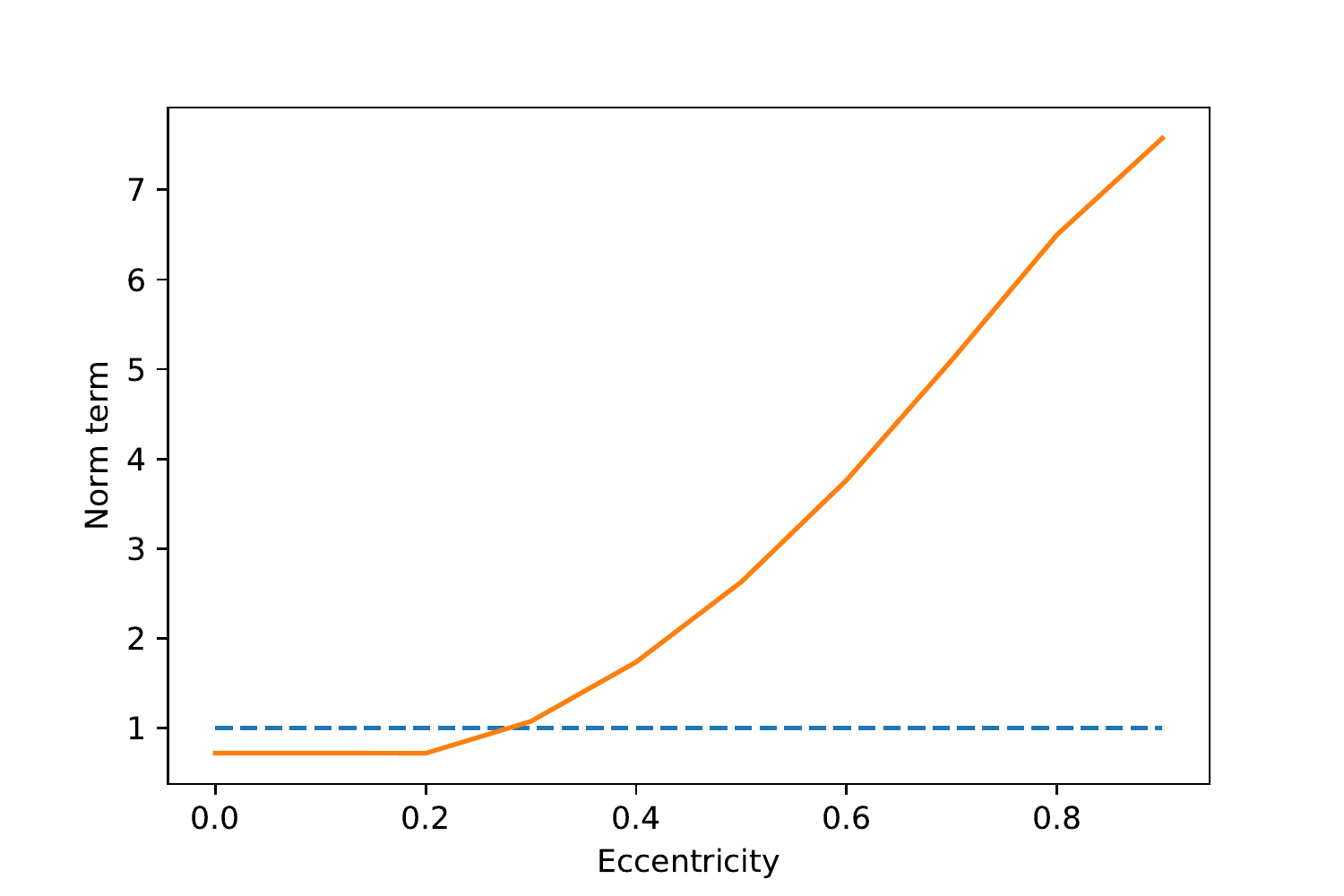}
\caption{Values of the norm terms for $\gamma = 3.5$ and
  eccentricities, $e \in [0,1)$, with a step size of $0.1$.}
\label{fig:all}
\end{subfigure}
\begin{subfigure}{0.5\textwidth}
\includegraphics[width=0.9\linewidth, height=5cm]{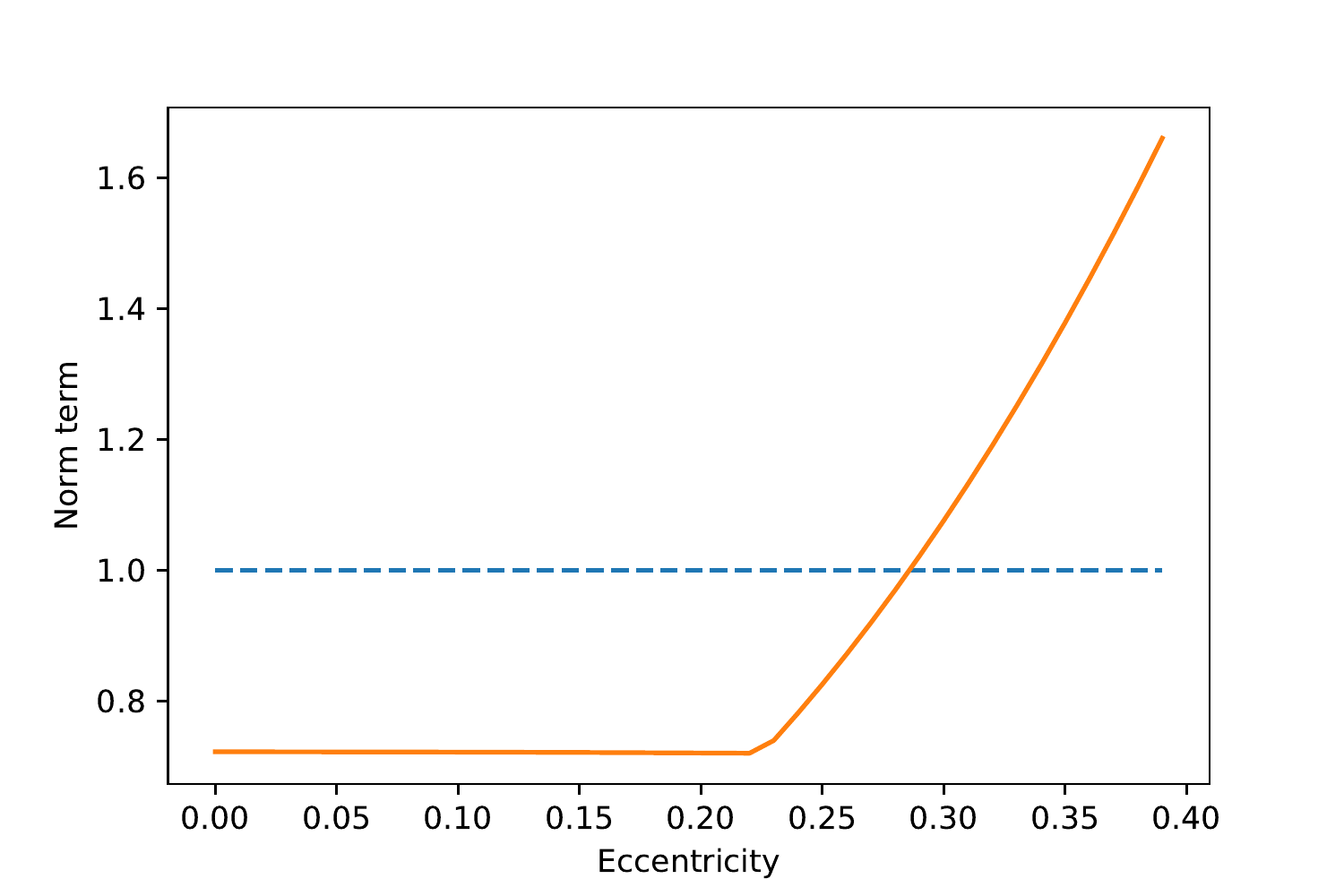}
\caption{Values of the norm terms for $\gamma = 3.5$ and eccentricities, $e \in [0,0.4)$,
  with a step size of 0.01.}
\label{fig:some}
\end{subfigure}

\caption{Plots for Norm terms vs eccentricities for the two cases the script was run. (a) The eccentricities, $e \in [0,1)$, with a step size of 0.1. The norm terms stay well below 1 until $e = 0.2$, but after 0.2 they start to increase and surpass 1 at $e = 0.3$, and continue to grow. (b) The eccentricities, $e \in [0,0.4)$, with a step size of 0.01. The norm terms stay well below 1 until $e = 0.2$, but after 0.2 they start to increase and surpass 1 at $e=0.28$, and continue to grow past $e = 0.3$.  }
\label{fig:both}
\end{figure}

While we found collision points for eccentricities, $0<e<1$ with step
size $0.01$, in the interest of computational time it took for the
script to run, we divided our analysis into stages. First, we found
the norm terms for $0<e<1$ with a step-size of $0.1$ (see
Figure~\ref{fig:all}), and successively refined the step-size in more
narrow ranges.

For the choice $\gamma = 3.5$, we observed that the norm terms ranged
from 0.722 to 7.573, as $e$ ranged from 0 to 0.9. The terms stayed
close to an approximate value of 0.72 for $e = 0, 0.1$ and $0.2$, but
then increase to 1.08 for $e = 0.3$ and 1.74 for $e=0.4$. (For the
norm terms generated for each eccentricity in this case, we refer the
reader to table \ref{tab:table_all}, Appendix A, Supplementary
Materials). In all cases for $e>0.22$, the largest norm term was reached
when $q = 3$ (for cases $e\leq 0.22$, the largest norm term was reached when $q=1$), after which they decayed (but appeared to increase again
after $q=20$. We cannot exclude that this is caused by numerical
instability\footnote{This may be due to the system's sensitivity to
  initial conditions. If the conditions were highly accurate, we would
  see the terms to decay continuously.} and so we stopped the
computation once the terms were under 0.5. In the case where the terms
never went below 1, such as the case was for $e=0.9$, the computation
was stopped at $q = 30$.

Figure \ref{fig:all} shows that after $e=0.4$, the terms are well past
1, and so we took a closer look at eccentricities below $0.4$. Hence,
we carried out the computations for $0<e<0.4$ with a step-size of
$0.01$ to better locate the $e$ at which the terms crosses the value
$1$ (see Figure~\ref{fig:some}).

We found that in this case, the terms stayed relatively close to a
value of $0.72$ (obtained for $q = 1$) until $e=0.22$, indicated by
the almost horizontal line in Figure~\ref{fig:some}, after which they
increased to $0.82$ for $e=0.25$, surpassing 1 at $e=0.29$. (For the
norm terms generated for each eccentricity in this case, we refer the
reader to table \ref{table_some}, Appendix A, Supplementary
Materials). The norm terms were observed to be equal to $0.97$ at
$e=0.28$ and equal to $1.022$ at $e=0.29$, indicating that the matrix
$\mathcal{T}_{q,j}$ might no longer be invertible once $e > 0.28$.  We
then modified the value of $\gamma$ and explored a narrower range of
eccentricities. A choice of $\gamma=3.1$ gave a maximum norm term of
approximately $0.85$ for this eccentricity. It also gave a norm term
of $0.968$ for $e=0.32$. We investigated what would happen if we were
to choose $\gamma=3.01$, but this only extended the range of
eccentricities to $e\le 0.33$ as it gave a norm term of approximately
$0.98$ (see Figure \ref{fig:different_gamma}). All of these
computations took a total of approximately 8 hours to run\footnote{The
  computation with $\gamma=3.1$ was run for $e\in[0.25,0.4]$, and with
  $\gamma=3.01$ was run for $e\in[0.32,0.4]$.}. Higher eccentricities
were not included in the computations involving smaller values for
$\gamma$, as computational time becomes an issue.


An interesting result we noticed
was that the maximum norm term was consistently achieved for $q=3$ for
$e>0.22$, and it was also observed that after $e \geq 0.15$, the terms tended to decrease for $q=2$
but then increase for $q=3$ before steadily decreasing again.

\subsection{Remarks and Future Suggestions}

In this section we make some remarks on our results described in the
previous section, as well as some recommendations that could improve
the results of this report.
\begin{enumerate}
  \item \label{item:minor-axis} The ellipse has two axis of symmetry.
  We performed a similar analysis considering perturbations preserving
  the symmetry along the minor axis and found no difference in the
  outcome: ellipses appears to be spectrally rigid at least up to
  $e = 0.3$
\item

We note here that our choice of $\gamma$ was rather arbitrary. The method
described in \cite{dsr} only requires $\gamma$ to be such that
$\gamma \in (3,4)$ for which $\norm{\tilde\cT_{q,j}-\Id} < 1$.  As
shown in~\cite{dsr}, the norm terms decay with a rate of
$Cq^{\gamma-3}$, hence the smaller the value for $\gamma$ used, the
slower the decay, and correspondingly, the longer the computational
time.   For the norm terms generated for each eccentricity in these
cases, we refer the reader to tables \ref{tab:table_3.1gamma} and
\ref{tab:table_3.01gamma}, Appendix A, Supplementary Materials.

\begin{figure}[ht]
\includegraphics[width=8cm]{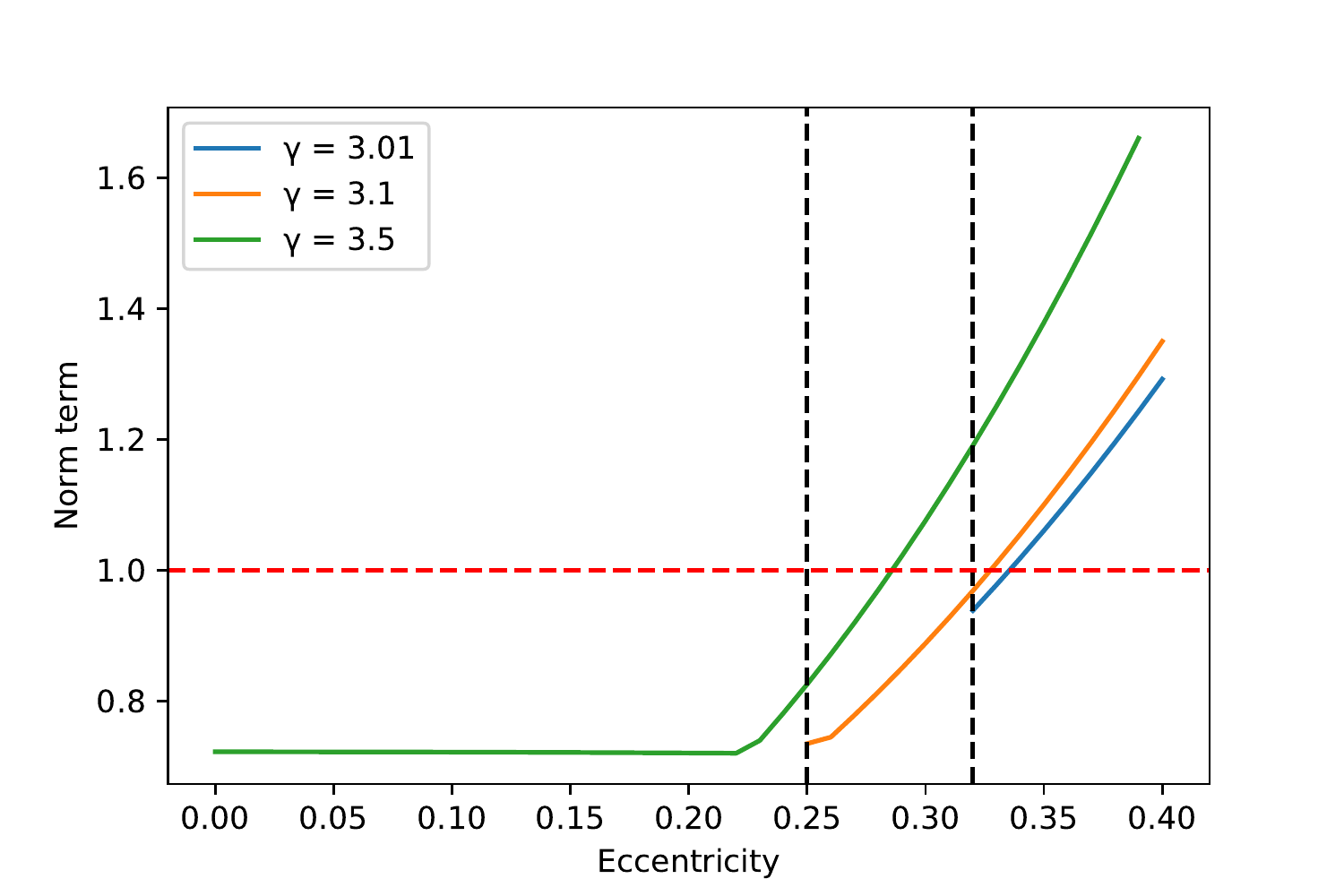}
\caption{Norm terms vs eccentricity with different values of $\gamma$: 3.5 (green), 3.1 (orange), 3.01 (blue). The dotted lines denote the eccentricity at which each was run, $\gamma=3.1$ was run for $0.25<e<0.4$ and $\gamma=3.01$ was run for $0.32<e<0.4$. The norm terms exceeded 1 for $e>0.32$ with $\gamma=3.1$ and for $e>0.33$ with $\gamma=3.01$.}
\label{fig:different_gamma}
\end{figure}


\item So far in our work we have chosen $S = \Id$ to determine
whether or not the operator was invertible, by means of
Proposition~\ref{p_injectivity}.  However, one could choose different
operators $S$; a natural choice would be $S = \cT^{0}$ (recall that
$\cT^{0}$ is the operator $\cT$ for the circle) . This should improve
the range of eccentricities for which we are able to provide numerical
evidence of spectral rigidity, but most likely the improvement would
not be substantial.

More interestingly, the finding that the highest norm term was
consistently achieved for $q=3$ could be further exploited to
construct an operator $S$ which could be more efficient in comparing
with $\cT$.  In our work, after the norm terms for $q > 3$ seemed to
decay quickly and the term $q = 3$ would often be the only term above
$1$.  Constructing a suitable operator would possibly lead to
substantial improvements in the range of confidence of our numerical
explorations.
\end{enumerate}
\newpage

\section{Appendix A: Supplementary Materials}

\begin{table}[ht]
    \centering
     \caption{Table showing the Norm terms generated for $e \in [0,1)$ with a step size of 0.1.}
    \begin{tabular}{|c|c|}
    \hline
     Eccentricity  & Norm Term  \\
     \hline
     0.00 & 0.7220 \\
     \hline
     0.10 & 0.7215 \\
     \hline
     0.20 & 0.7202 \\
     \hline
     0.30 & 1.0757 \\
     \hline
     0.40 & 1.7370 \\
     \hline
     0.50 & 2.6304 \\
     \hline
     0.60 & 3.7642 \\
     \hline
     0.70 & 5.1015 \\
     \hline
     0.80 & 6.4986 \\
     \hline
     0.90   & 7.5732 \\
    \hline
    \end{tabular}
    \label{tab:table_all}
\end{table}

\begin{longtable}[c]{|c|c|}
\caption{Table showing the norm terms generated for $e \in [0,0.4)$ with a step size of 0.01.\label{table_some}}\\
    \hline
       Eccentricity  & Norm Term  \\
       \hline
       0.00 &	0.7220 \\
       \hline
      0.01 &	0.7220 \\
      \hline
0.02 &	0.7220 \\
\hline
0.03 &	0.7219 \\
\hline
0.04 &	0.7219 \\
\hline
0.05 &	0.7219 \\
\hline
0.06 &	0.7218 \\
\hline
0.07 &	0.7218 \\
\hline
0.08 &	0.7217 \\
\hline
0.09 &	0.7216 \\
\hline
0.10 &	0.7215 \\
\endfirsthead
\hline
0.11 &	0.7214 \\
\hline
0.12 &	0.7213 \\
\hline
0.13 &	0.7212 \\
\hline
0.14 &	0.7211 \\
\hline
0.15 &	0.7210 \\
\hline
0.16 &	0.7208 \\
\hline
0.17 &	0.7207 \\
\hline
0.18 &	0.7205 \\
\hline
0.19 &	0.7204 \\
\hline
0.20 &	0.7202 \\
\hline
0.21 &	0.7200 \\
\hline
0.22 &	0.7198 \\
\hline
0.23 &	0.7393 \\
\hline
0.24 &	0.7814 \\
\hline
0.25 &	0.8254 \\
\hline
0.26 &	0.8714 \\
\hline
0.27 &	0.9194 \\
\hline
0.28 &	0.9695 \\
\hline
0.29 &	1.0216 \\
\hline
0.30 &	1.0757 \\
\hline
0.31 &	1.1320 \\
\hline
0.32 &	1.1904 \\
\hline
0.33 &	1.2510 \\
\hline
0.34 &	1.3137 \\
\hline
0.35 &	1.3786 \\
\hline
0.36 &	1.4458 \\
\hline
0.37 &	1.5151 \\
\hline
0.38 &	1.5868 \\
\hline
0.39 &	1.6607 \\
\hline
\end{longtable}

\begin{table}[ht]
    \centering
     \caption{Table showing the Norm terms generated for $e \in [0.25,0.4]$ with $\gamma=3.1$.}
    \begin{tabular}{|c|c|}
    \hline
     Eccentricity  & Norm Term  \\
     \hline
     0.25 & 0.7345 \\
     \hline
     0.26 & 0.7444 \\
     \hline
     0.27 & 0.7782 \\
     \hline
     0.28 & 0.8133 \\
     \hline
     0.29 & 0.8499 \\
     \hline
     0.30 & 0.8879 \\
     \hline
     0.31 & 0.9274 \\
     \hline
     0.32 & 0.9683 \\
     \hline
     0.33 & 1.0107 \\
     \hline
     0.34 & 1.0546 \\
     \hline
     0.35 & 1.1000 \\
     \hline
     0.36 & 1.1469 \\
     \hline
     0.37 & 1.1954 \\
     \hline
     0.38 & 1.2453 \\
     \hline
     0.39 & 1.2969 \\
     \hline
     0.40 & 1.3500 \\
    \hline
    \end{tabular}
    \label{tab:table_3.1gamma}
\end{table}
\vspace{2cm}
\newpage
\begin{table}[ht]
    \centering
     \caption{Table showing the Norm terms generated for $e \in [0.3,0.4]$ with $\gamma=3.01$.}
    \begin{tabular}{|c|c|}
    \hline
     Eccentricity  & Norm Term  \\
     \hline
     0.32 & 0.9382 \\
     \hline
     0.33 & 0.9775 \\
     \hline
     0.34 & 1.0183 \\
     \hline
     0.35 & 1.0604 \\
     \hline
     0.36 & 1.1039 \\
     \hline
     0.37 & 1.1488 \\
     \hline
     0.38 & 1.1951 \\
     \hline
     0.39 & 1.2429 \\
     \hline
     0.40 & 1.2921 \\
    \hline
    \end{tabular}
    \label{tab:table_3.01gamma}
\end{table}

Link to the code:
\href{https://github.com/shanzaayub/Spectral-Rigidity-Ellipse}{GitHub repository}

\end{document}